\documentclass[12pt]{article}
\usepackage{amsmath,amsthm,amssymb, amstext, amsopn, amsxtra}
\begin{document}
\newcommand{\nt}{\noindent}
\newcommand{\bs}{\bigskip}
\newcommand{\ms}{\medskip}
\newcommand{\mk}{\medskip}
\newcommand{\sk}{\smallskip}

\[    \]
\begin{center}
\Large{\bf Finite Groups Embeddable in Division Rings}

\bigskip
\large{\bf T. Y. Lam } \\
\end{center}

\bigskip
\begin{center}
{\bf \S1. Introduction}
\end{center}

\mk
Finite groups that are embeddable in the multiplicative groups of division 
rings $\,K\,$ were completely determined by S.~A.~Amitsur [Am] in 1955.  
In case $\,K\,$ has characteristic $\,p>0$, the only possible finite 
subgroups of $\,K^*\,$ are cyclic groups, according to a theorem of 
I.~N.~Herstein [He].  Thus, the only interesting case is when $\,K\,$ has 
characteristic $\,0$; that is, when $\,K\supseteq {\Bbb Q}\,$.

\bs
In [He], Herstein conjectured that odd-order subgroups of division rings
$\,K\,$ were cyclic, and he proved this to be the case when $\,K\,$ is the
division ring of the real quaternions.  Herstein's conjecture was settled
negatively in [Am].  As part of his complete classification of finite
groups in division rings, Amitsur showed that the smallest noncyclic
odd-order group that can be embedded in a division ring is one of order 
$63$ (and this group is unique).

\bs
Amitsur's paper is daunting to read as it is long and technically 
complicated.  In lecturing to a graduate class on division rings, I tried 
to find a simple reason for the ``first exceptional odd order'' $63$ (to
Herstein's conjecture).
After some work, I did come up with a reason that was simple enough to be 
explained to my class, without having to go through {\it any part\/} of 
Amitsur's paper.  Furthermore, the method I used led easily 
to the {\it second\/} exceptional odd order, $117$ (which was not mentioned
in Amitsur's paper).  Since this line of reasoning did not seem to have 
appeared in the literature before, I record it in this short note.
To better motivate the results discussed here, I have also included
a quick exposition on the beginning part of the theory of finite subgroups 
of division rings. As a result, this paper can be read independently of [Am].

\sk

\bs\nt

\begin{center}
{\bf \S2. Background on Amitsur Groups}
\end{center}

\mk
Let us say that a finite group $\,G\,$ is {\it Amitsur\/} if $\,G\,$ can be 
embedded in the multiplicative group $\,K^*\,$ of some division ring $\,K$.
We start with some background information on such Amitsur groups.  First 
the following easy observation.

\bs\nt
{\bf (2.1) Theorem.} {\it Any abelian Amitsur group $\,G\,$ is cyclic.}

\bs\nt
{\bf Proof.} Say $\,G\subseteq K^*$, where $\,K\,$ is a division ring.
If $\,F\,$ is the center of $\,K$, then the $\,F$-span $\,L\,$ of $\,G\,$ 
is a finite-dimensional $\,F$-algebra.  Since $\,L\,$ is a commutative
domain, it must be a field.  From field theory (see, e.g., [Ja\,:\,p.\,132]),
it then follows that $\,G\subseteq L^*\,$ is cyclic.
\ \ \ \ \ \ \ \ {\bf QED}

\bs
Next, we classify Amitsur $\,p$-groups.  (Here, $\,p\,$ denotes any 
prime number, unrelated to the characteristic of our division rings.)
This is possible thanks to some standard results in $\,p$-group theory.

\bs\nt
{\bf (2.2) Theorem.} {\it Let $\,G\,$ be an Amitsur $\,p$-group. Then
$\,G\,$ is cyclic if $\,p\,$ is odd, and $\,G\,$ is cyclic or 
dicyclic if $\,p=2$.}

\bs
Recall that a group is {\it dicyclic\/} (of order $\,4n$) if it is generated 
by two elements $\,x,\,y\,$ subject to the relations $\,x^{2n}=1,\;y^2=x^n$, 
and $\,yxy^{-1}=x^{-1}$, where $\,n\geq 2$.  Dicyclic groups are
also known as {\it generalized quaternion groups\/}; they are $\,2$-groups 
iff $\,n\,$ is a power of $\,2$.  It is easy to see that any dicyclic group 
$\,\langle x,\,y \rangle\,$ as above is embeddable in the division ring 
of real quaternions, so the conclusion in (2.2) is indeed the best possible.

\bs\nt
{\bf Proof of (2.2).} In the theory of $\,p$-groups, there is a well-known 
result which guarantees that, if a $\,p$-group has a unique subgroup of 
order $\,p$, then it satisfies the conclusion of the present theorem; see
Th.~12.5.2 in [Ha].  Thus, it suffices for us to verify that, if 
$\,G\neq \{1\}$, $\,G\,$ has a unique subgroup of order $\,p$.  Say
$\,G\subseteq K^*$, where $\,K\,$ is a division ring, with center $\,F$.

\sk
By elementary group theory, we know that $\,G\,$ has a nontrivial 
center.  Fix a central element $\,x\in G\,$ of order $\,p$, so
$\,|\langle x\rangle |=p$.  If $\,\langle y\rangle \subseteq G\,$ is
any subgroup of order $\,p$, then $\,L:=F(x,y)\,$ is a field.  In this
field, the equation $\,t^p-1=0\,$ has at most $\,p\,$ solutions, so
we must have $\,\langle y\rangle = \langle x\rangle $, as desired. 
\ \ \ \ \ \ \ \ {\bf QED}

\bs
In order to get some general results for non $\,p$-groups, we introduce
the following class of finite groups, where $\,m,\,n,\,r\,$ are natural
numbers:
$$ G_{m,n,r}=\langle a,\,b\,: \;\;a^m=b^n=1,\;\;bab^{-1}=a^r \rangle \,.
\leqno (2.3) $$
Whenever we use this notation, it will be assumed that $\,r^n\equiv 1\;
(\mbox{mod}\;m)$.  This assumption guarantees that $\,G_{m,n,r}\,$ is a 
semidirect direct product of $\,\langle a\rangle\,$ and 
$\,\langle b\rangle\,$, with $\,|\langle a\rangle|=m\,$ and
 $\,|\langle b\rangle|=n$.  In particular, $\,|G_{m,n,r}|=mn$.
If necessary, we may assume $\,0\leq r<m$; but this is not always essential.

\bs
Let us say that a finite group $\,G\,$ is {\it Sylow-cyclic\/} if all of its Sylow 
subgroups are cyclic.\footnote{Such groups are also called {\it $\,Z$-groups\/} 
in the literature, after Hans Zassenhaus. We use here the term 
``Sylow-cyclic'' since it is completely self-explanatory.}  
For instance, any $\,G\,$ with a square-free 
order is Sylow-cyclic, since each of its Sylow groups has a prime order. The
following classical result of H\"older, Burnside, and Zassenhaus (see Th.~9.4.3 
in [Ha]) gives a complete determination of all Sylow-cyclic finite groups.

\bs\nt
{\bf (2.4) H\"older-Burnside-Zassenhaus Theorem.} {\it A finite group $\,G\,$ is
Sylow-cyclic iff $\,G\cong G_{m,n,r}\,$ for some $\,m,\,n,\,r\,$ with
$\,(m,\,n(r-1))=1$. (In particular, Sylow-cyclic groups are solvable.)}

\bs
Coupling this powerful classification result with (2.2), we deduce the following 
important information on a class of Amitsur groups, including all the ones of
odd order.

\bs\nt
{\bf (2.5) Theorem.} {\it If $\,G\,$ is an Amitsur group with order not divisible
by $8$, then $\,G\,$ is Sylow-cyclic, and hence $\,G\cong G_{m,n,r}\,$ for some
$\,m,\,n,\,r\,$ with $\,(m,\,n(r-1))=1.$}

\bs
According to this result, Amitsur groups with order not divisible by $8$ are
``fairly close'' to being cyclic: they are at least {\it metacyclic\/}, i.e.,
extensions of cyclic groups by cyclic groups.  To find all odd-order Amitsur
groups, for instance, we must determine, for $\,m,\,n\,$ odd, which of the
metacyclic groups $\,G_{m,n,r}\,$ (with $\,(m,\,n(r-1))=1\,$) are Amitsur.
In the next section, we shall present some partial results in this direction,
not only for odd-order groups, but for the groups $\,G_{m,n,r}\,$ in (2.3)
in general.  Note that the groups  $\,G_{m,n,r}\,$ in this paper are ``simpler 
versions'' of the groups $\,G_{m,r}\,$ studied by Amitsur, in that the
$\,G_{m,n,r}\,$'s are given as semidirect products, while Amitsur's 
$\,G_{m,r}\,$'s are not.  As we shall see in the next section, the semidirect 
product representation of $\,G_{m,n,r}\,$ makes it possible for us to obtain 
quickly some necessary conditions for such a group to be Amitsur.

\sk

\bs
\begin{center}
{\bf \S3. Necessary Conditions for the Embeddability of $\,G_{m,n,r}$}
\end{center}

\mk
In this section, we shall study the groups $\,G_{m.n,r}\,$ that are embeddable
in division rings.  (We do not need to assume the condition
$\,(m,\,n(r-1))=1$.)  The main result here is first stated in the following
negative fashion.

\bs\nt
{\bf (3.1) Theorem.} {\it Let $\,G=G_{m,n,r}\,$, where $\,n>1\,$ and}
$\,o\,(\overline{r})\,$ ({\it the order of $\,\overline{r}\,$ in the unit group} 
$\,\mbox{U}({\Bbb Z}_m)\,$) {\it is exactly $\,n$.  Then $\,G\,$ is not Amitsur.}

\bs\nt
{\bf Proof.} Assume, instead, that $\,G\subseteq K^*$, where $\,K\,$ is
a division ring.  We use the presentation (2.3) for $\,G$. From
$\,bab^{-1}=a^r$, we get $\,b^iab^{-i}=a^{r^i}$.  Thus, $\,b^ia=a^{r^i}b$,
and $\,a^{-1}b^ia=a^{r^i-1}b^i$.  From $\,b^n=1$, we get
$\,b^{n-1}+\cdots+b+1=0\,$ in $\,K$.  Conjugating this equation by $\,a$,
we get
$$ a^{r^{n-1}-1}b^{n-1}+a^{r^{n-2}-1}b^{n-2}+\cdots +a^{r-1}b+1=0. \leqno (3.2)$$
Subtracting the last two equations and cancelling $\,b\,$ (from the right)
gives
$$ \bigl(a^{r^{n-1}-1}-1\bigr)b^{n-2}+\bigl(a^{r^{n-2}-1}-1\bigr)b^{n-3}+
\cdots +(a^{r-1}-1)=0. \leqno (3.3)$$
Conjugating this by $\,a\,$ again and subtracting (3.3) from the resulting
equation, we get after another cancellation of $\,b$:
$$ \bigl(a^{r^{n-1}-1}-1\bigr)\bigl(a^{r^{n-2}-1}-1\bigr)b^{n-3}+\cdots
+(a^{r^2-1}-1)(a^{r-1}-1)=0.$$
Carrying this process to the bitter end produces the equation
$$ \bigl(a^{r^{n-1}-1}-1\bigr)\bigl(a^{r^{n-2}-1}-1\bigr)\cdots
(a^{r^2-1}-1)(a^{r-1}-1)=0 \in K.  \leqno (3.4)$$
Therefore, $\,a^{r^i-1}=1\,$ for some $\,i\in [1,\,n-1]$.  But then
$\,r^i\equiv 1\,(\mbox{mod}\;m)$, which contradicts $\,o\,(\overline{r})=n$.
\,\,\,\,\,{\bf QED}

\bs 
In my lectures, I referred to the above as the ``triple $(-1)$-proof'',
because of the three remarkable layers of $(-1)\,$'s in the first factor
of the key equation (3.4)!  Note that this proof actually showed that, under the 
assumption that $\,o\,(\overline{r})=n>1$, $\,G_{m,n,r}\,$ cannot be 
multiplicatively embedded in any domain (that is, a nonzero ring without 
0-divisors).  This, however, is not a real improvement, since it is easily seen 
that, once a finite (multiplicative) group $\,G\,$ is embeddable in a domain, 
then $\,G\,$ is in fact embeddable in a division ring.

\bs 
We now record some consequences of (3.1).

\bs\nt
{\bf (3.5) Corollary.} {\it If $\,G=G_{m,n,r}\,$} ({\it generated by $\,a,\,b\,$ 
as in} $(2.3)$) {\it is Amitsur, then the largest square-free subgroup of
$\,\langle b\rangle\,$ must be central in $\,G$.} {\it In particular, if
$\,G=G_{m,n,r}\,$ is Amitsur and $\,n\,$ is square-free, then $\,G\,$ is cyclic.}

\bs\nt
{\bf Proof.} For any prime $\,\ell\,$ dividing $\,n=o(b)$, let 
$\,\langle c \rangle\,$ be the $\,\ell$-Sylow group of $\,\langle b \rangle\,$, 
say of order $\,\ell ^t$.  Then
$\,cac^{-1}=a^s\,$ for some $\,s$, and 
$\,\langle a\rangle  \langle c\rangle \cong G_{m,\ell ^t,s}\,$ is still Amitsur.  
By (3.1), we must have
$\,o\,(\overline{s})< \ell ^t\,$ in $\,\mbox{U}({\Bbb Z}_m)$.  This means that
$\,c^{\ell ^{t-1}}\,$ acts trivially on $\,a\,$ (by conjugation), and hence
$\,c^{\ell ^{t-1}}\,$ is central in $\,G$. From this, the first conclusion in the
Corollary follows.  If $\,n\,$ happens to be square-free, this conclusion means 
that $\,b\,$ is central in $\,G$.  Thus, $\,G\,$ is abelian, and hence cyclic
by (2.1). \ \ \ \ \ \ \ \ \ {\bf QED}

\bs\nt
{\bf (3.6) Corollary.} {\it Let $\,G\,$ be an Amitsur group of order
$\,p^{t}q_1\cdots q_k\,$, where $\,k\geq 1$, and $\,p>q_1>\cdots >q_k\,$ are primes.
Then $\,G\,$ is cyclic.}

\bs\nt
{\bf Proof.} Since $\,k\geq 1$, $\,p\,$ is an odd prime. Thus, $\,G\,$ is 
Sylow-cyclic, and hence solvable by (2.4).  We proceed by induction on $\,k$.
Note that the solvability of $\,G\,$ implies, by Philip Hall's Theorem, that it 
has a subgroup $\,H\,$ of order $\,m:=p^tq_1\cdots q_{k-1}$ (see [Ha:~Th.~9.3.1]).
This Hall subgroup is still Amitsur, and hence cyclic by the inductive hypothesis. 
Also, since $\,[G:H]=q_k\,$ is the smallest prime divisor of $\,|G|$, $\,H\,$
must be normal in $\,G$.  Thus, by taking a generator of $\,H\,$ and an element 
of order $\,q_k\,$ in $\,G$, we can represent $\,G\,$ as $\,G_{m,q_k,r}\,$ for 
some $\,r$.  Since $\,q_k\,$ is a prime, (3.5) implies that $\,G\,$ is cyclic.  
Note that this argument is also sufficient to get the induction started, so the
proof is complete. \ \ \ \ \ \ \ \ {\bf QED}

\bs
A special case of (3.6) is that {\it any square-free Amitsur group is cyclic.}
This recaptures Corollary 5 on p.\,384 of [Am].

\bs
Note that, in (3.1) and (3.5), the integers $\,m,\,n\,$ and $\,r\,$ were
arbitrary.  To get some explicit numerical results, let us specialize these 
results to odd-order groups.

\bs\nt
{\bf (3.7) Theorem.} {\it If $\,G\,$ is an Amitsur group of odd order $< 171$,
then $\,G\,$ is cyclic except possibly when} $\,|G|=63\;\,or \;\,117.$

\bs\nt
{\bf Proof.} If $\,|G|\,$ is either a prime power, or square-free, or of 
the form $\,p^{t}q\,$ where $\,p>q\,$ are primes, the foregoing results 
imply that $\,G\,$ is cyclic.  Among odd integers from $1$ to $169$,
the only ones {\it not\/} of any of the above forms are
$$ 5\!\cdot\! 3^2=45, \;\;\, 7\!\cdot\! 3^2=63, \;\;\, 11\!\cdot \!3^2=99,
\;\;\, 13\!\cdot\! 3^2=117, \;\;\, 5\!\cdot\! 3^3=135 \;\;\, \mbox{and}
\;\;\, 17\!\cdot\! 3^2=153.$$
If $|G|=45$ or $99$, $\,G\,$ is easily seen to be abelian (by Sylow
theory), and hence cyclic (by (2.1)).  Next, consider $\,|G|=135$. Since
$\,G\,$ is Amitsur, we can represent it (thanks to (2.5)) in the form 
$\,G_{m,n,r}\,$ with $\,(m,\,n)=1$.  We may assume $\,m,\,n>1$, so we have 
$\,m=5\,$ or $\,3^3$.  But $\,\mbox{U}({\Bbb Z}_5)\,$ has no element of order 
$\,3$, and $\,\mbox{U}({\Bbb Z}_{3^3})\,$ has no element of order $\,5$.  
Therefore, $\,G\,$ must be abelian and hence cyclic.  The case $\,|G|=153\,$
can be treated similarly  (with exactly the same conclusion).  Thus, $\,G\,$ 
is cyclic in all cases, except possibly when $\,|G|=63\,$ or $\,117$. 
\ \ \ \ \ \ \ \ \ {\bf QED}

\bs
What about $\,|G|=63$?  There are, up to isomorphism, four groups of order 63;
they are:
$$ {\Bbb Z}_{63},\;\;\;\;\; {\Bbb Z}_3 \times {\Bbb Z}_{21}, \;\;\;\;\;
{\Bbb Z}_3\times G_{7,3,2}\,, \;\;\;\;\mbox{and} \;\;\;\; G_{7,9,2}.$$
The first group is cyclic, and hence Amitsur.  The second and third ones are
not Sylow-cyclic, and hence not Amitsur.  The remaining question is whether the
non-cyclic (but Sylow-cyclic) $\,G_{7,9,2}\,$ is Amitsur.  Note that, since
$\,o\,(\overline{2})=3<9\,$ in $\,\mbox{U}({\Bbb Z}_7)$, (3.1) does not apply
to this group.  Thus, there {\it is\/} a chance that $\,G_{7,9,2}\,$ is Amitsur.
Indeed, the embeddability of this group was proved by Amitsur, and subsequent
expositions on this have been given by C.~Ford [Fo] and J.~Dauns [Da].  For 
the sake of completeness, we give a quick sketch for the embeddability of
$\,G_{7,9,2}\,$ below.

\bs
Let $\,L={\Bbb Q}\,(\zeta)\,$ where $\,\zeta\,$ is a primitive 21st root of unity,
and let $\,\sigma\in \mbox{Gal}\,(L/{\Bbb Q}\,)\,$ be defined by $\,\sigma(\zeta)
=\zeta^{16}$.  Then $\,o\,(\sigma)=3$, so $\,L\,$ is a cubic extension of the
fixed field $\,F:=L^{\sigma}$.  Note that $\,F\,$ contains the primitive
cubic root of unity $\,\omega:=\zeta^7$, since 
$$\,\sigma(\omega)=\zeta ^{7\cdot 16}=\zeta^7=\omega.$$  
On the other hand, for the primitive 7th root of unity $\,a:=\zeta^3$, we have
$$\,\sigma(a)=\zeta^{3\cdot 16}=\zeta^{6}=a^2.$$ 
Now introduce a new symbol $\,b$, and form the cyclic $\,F$-algebra
$$ K:= L\oplus L\,b \oplus L\,b^2, \;\;\;\; \mbox{with}\;\;b^3=\omega, 
\;\;\;\mbox{and} \;\;\; b\,\ell=\sigma(\ell)\,b \;\;\; (\forall \;\,\ell \in L),$$
which has center $\,F$.  We have $\,\mbox{dim}\,_{\Bbb Q}\,F=\varphi(21)/3=4$, so 
$\,\mbox{dim}\,_{\Bbb Q}\,K=3^2\!\cdot\!4=36.$  In the algebra $\,K$, the elements
$\,a,\,b\,$  satisfy the equations
$$ a^7=1,\;\;\; b^3=\omega,\;\;\;b^9=\omega^3=1, \;\;\;\mbox{and}\;\;\;
ba=\sigma(a)b=a^2b.  $$
Therefore, $\,a,\,$ and $\,b\,$ generate a group $\,G\cong G_{7,9,2}\,$
in the group of units $\,\mbox{U}(K)$. {\it It can be checked that $\,K\,$ is a 
division $\,F$-algebra} (see, for instance, [Fo]), so this establishes the 
noncyclic $\, G_{7,9,2}\,$ as an Amitsur group (of order $63$).  The center 
$\,F\,$ of the division ring $\,K\,$ is easily seen to be the biquadratic 
extension $\,{\Bbb Q}\bigl(\sqrt{-3},\,\sqrt{-7}\,\bigr)\,$ of the rationals.

\bs
As for the case $\,|G|=117$, we have a noncyclic candidate $\,G_{13,9,9}\,$
(noting that $\,o\,(\overline{9})=3\,$ in $\,\mbox{U}({\Bbb Z}_{13})$).  This
group can be shown to be Amitsur in the same way as $\,G_{7,9,2}\,$ was.

\bs
A substantial part of the argument in showing that a certain group is Amitsur
is to establish that a certain $\,{\Bbb Q}\,$-algebra is a division algebra.
This is by no means routine.  In Amitsur's paper, this part of the work is 
handled by number-theoretic tools, such as the Hasse Norm Theorem and the 
Hasse-Brauer-Noether Theorem.  Thus, to fully understand Amitsur's work would 
require considerable preparation in algebraic number theory.

\bs
In closing, we should mention that it is also of interest to {\it fix\/} a field 
$\,F\,$ (of characteristic zero), and to determine the finite groups $\,G\,$ 
that are embeddable in $\,K^*\,$ for some finite-dimensional {\it central\/} 
$\,F$-division algebra $\,K$. This problem has been studied by B.~Fein and 
M.~Schacher [FS], who called such $\,G\,$ ``$F$-adequate groups''.  The embedding 
constructed for $\,G=G_{7,9,2}\,$ above showed that $\,G\,$ is
$\,{\Bbb Q}\bigl(\sqrt{-3},\,\sqrt{-7}\,\bigr)$-adequate, but Fein and Schacher
showed that $\,G\,$ is already $\,{\Bbb Q}\bigl(\sqrt{-3}\,\bigr)$-adequate.
This is remarkable, since Fein and Schacher have also proved the following
interesting positive result on Herstein's conjecture: 
{\it If} $\,[F:{\Bbb Q}\,]\leq 2$, {\it then any $\,F$-adequate odd-order group 
is cyclic, except when} $\,F={\Bbb Q}\bigl(\sqrt{-3}\,\bigr)\,$!

\mk\nt
University of California \\
Berkeley, Ca 94720

\sk\nt
lam@math.berkeley.edu

\end{document}